\documentclass[a4]{article}

\font\ess=eufb10 scaled\magstep3 

\font\es=eufm10

\def\gC{\mbox{\es {C}}}

\def\gJ{\mbox{\es {J}}}

\def\tr{\mbox{\rm {tr}}}

\def\diag{\mbox{\rm {diag}}}
\def\Iso{\mbox{\rm {Iso}}}

\def\ov{\overline}

\def\dfrac#1#2{\displaystyle \frac{#1}{#2}}

\def\C{\mbox{\boldmath $C$}}

\def\H{\mbox{\boldmath $H$}}

\def\R{\mbox{\boldmath $R$}}

\def\0{\mbox{\boldmath {0}}}    
\def\1{\mbox{\boldmath {1}}}      
\def\2{\mbox{\boldmath {2}}}      
\def\3{\mbox{\boldmath {3}}}      
\def\4{\mbox{\boldmath {4}}}      
\def\5{\mbox{\boldmath {5}}}      
\def\6{\mbox{\boldmath {6}}}      
\def\7{\mbox{\boldmath {7}}}      
\def\8{\mbox{\boldmath {8}}}      
\def\9{\mbox{\boldmath {9}}}      
\def\a{\mbox{\boldmath $a$}}
\def\b{\mbox{\boldmath $b$}}

\def\lgJ{\mbox{\ess {J}}}

\begin{document}
\baselineskip=14pt

\begin{center}
\Large{\bf Constructive diagonalization of an element $X$ of the Jordan algebra $\lgJ$ by the exceptional group $F_4$}
\end{center}
\vspace{3mm}

\vspace{3mm}

\begin{center}
     \large{Takashi Miyasaka and Ichiro Yokota}
\end{center}
\vspace{4mm}

I. We know that any element $X$ of the exceptional Jordan algebra $\gJ$ is transformed to a diagonal form by the compact exceptional Lie group $F_4$. However, its proof is used the method which is reduced a contradiction. In this paper, we give a direct and constructive proof. 
\vspace{2mm}

Let $\H$ be the field of quaternions and $\gC = \H \oplus \H e_4, {e_4}^2 = -1$  the division Cayley algebra. For $K = \H, \gC$, let $\gJ(3, K) = \{X \in M(3, K) \, | \, X^* = X \}$ with the Jordan multiplicaton $X \circ Y$, the inner product $(X, Y)$ and the Freudenthal multiplication $X \times Y$ espectively by  
$$
\begin{array}{c}
      X \circ Y = \dfrac{1}{2}(XY + YX), \quad (X, Y) = \tr(X \circ Y),
\vspace{1mm}\\
   X \times Y = \dfrac{1}{2}(2X \circ Y - \tr(X)Y - \tr(Y)X + (\tr(X)\tr(Y) - (X, Y))E) 
\end{array}$$
(where $E$ is the $3 \times 3$ unit matrix).
\vspace{2mm}

The simply connected compact Lie group $F_4$ is defined by
\begin{eqnarray*}
   F_4 \!\!\! &=& \!\!\! \{ \alpha \in \Iso_R(\gJ(3, \gC)) \, | \, \alpha(X \circ Y) = \alpha X \circ \alpha Y \}
\vspace{1mm}\\
       \!\!\! &=& \!\!\! \{ \alpha \in \Iso_R(\gJ(3, \gC)) \, | \, \alpha(X \times Y) = \alpha X \times \alpha Y \}. 
\end{eqnarray*}
Then, we have the following Theorem ([1],[2]).
\vspace{3mm}

{\bf Theorem 1.} {\it Any element $X$ of $\gJ(3, \gC)$ can be transformed to a diagonal form by some element} $\alpha \in F_4$:
$$
     \alpha X = \pmatrix{\xi_1 & 0 & 0 \cr
                         0 & \xi_2 & 0 \cr
                         0 & 0 & \xi_3}, \quad \xi_i \in \R. $$
To give a constructive proof of this theorem 1, we will prepare some elements of $F_4$.

(1) Let $Sp(3) = \{A \in M(3, \H) \, | \, A^*A = E \}$. We shall show that the group $F_4$ contains $Sp(3)$ as subgroup : $Sp(3) \subset F_4$. An element $X = \gJ(3, \gC)$ is expressed by
$$
   X = \pmatrix{\xi_1 & x_3 & \ov{x}_2 \cr
                            \ov{x}_3 & \xi_2 & x_1 \cr
                            x_2 & \ov{x}_1 & \xi_3}
     = \pmatrix{\xi_1 & m_3 & \ov{m}_2 \cr
              \ov{m}_3 & \xi_2 & m_1 \cr
              m_2 & \ov{m}_1 & \xi_3} 
     + \pmatrix{0 & a_3e_4 & -a_2e_4 \cr
                -a_3e_4 & 0 & a_1e_4 \cr                                                         a_2e_4 & -a_1e_4 & 0}, $$
where $x_i = m_i + a_ie_4 \in \H \oplus \H e_4 = \gC$. To such $X$, we associate an element 
$$
   \pmatrix{\xi_1 & m_3 & \ov{m}_2 \cr
              \ov{m}_3 & \xi_2 & m_1 \cr
              m_2 & \ov{m}_1 & \xi_3}  + (a_1, a_2, a_3) $$    
of $\gJ(3, \H) \oplus \H^3$. In $\gJ(3, \H) \oplus \H^3$, we define a multiplication $\times$ by
$$
   (M + \a) \times (N + \b) = \Big(M \times N - \dfrac{1}{2}(\a^*\b + \b^*\a)\Big) - \dfrac{1}{2}(\a N + \b M). $$
Since this multilication corresponds to the Freudenthal multiplication in $\gJ(3, \gC)$, hereafter, we identify $\gJ(3, \H) \oplus \H^3$ with $\gJ(3, \gC)$.
\vspace{2mm}

Now, we define a map $\varphi : Sp(3) \to F_4$ by
$$
   \varphi(A)(M + \a) = AMA^* + \a A^*, \quad M + \a \in \gJ(3, \H) \oplus \H^3 = \gJ(3, \gC).$$
It is not difficult to see that $\varphi$ is well-defined : $\varphi(A) \in F_4$. Since $\varphi$ is injective, we identify $Sp(3)$ with $\varphi(Sp(3))$ : $Sp(3) \subset F_4$ ([3]). 
\vspace{1mm}

(2) Let $G_2 = \{\alpha \in \Iso_R(\gC) \, | \, \alpha(xy) = (\alpha x)(\alpha y) \}$. The group $F_4$ contains $G_2$ as subgroup by the following way. We define a map $\phi : G_2 \to F_4$ by
$$
   \phi(\alpha)\pmatrix{\xi_1 & x_3 & \ov{x}_2 \cr
                        \ov{x}_3 & \xi_2 & x_1 \cr
                        x_2 & \ov{x}_1 & \xi_3}
   = \pmatrix{\xi_1 & \alpha x_3 & \ov{\alpha x_2} \cr
              \ov{\alpha x_3} & \xi_2 & \alpha x_1 \cr
                       \alpha x_2 & \ov{\alpha x_1} & \xi_3}. $$
It is not difficult to see that $\phi$ is well-defined : $\phi(\alpha) \in F_4$. Since $\phi$ is injective, we identify $G_2$ with $\phi(G_2)$ : $G_2 \subset F_4$. 

(3) For $D = \diag(a, \ov{a}, 1), a \in \gC, |a| = 1$, we define an $\R$-linear transformation $\delta_a$ of $\gJ(3, \gC)$ by
$$
   \delta_aX = D_aX\ov{D}_a = \pmatrix{\xi_1 & ax_3a & a\ov{x}_2 \cr
                            \ov{a}\ov{x}_3\ov{a} & \xi_2 & \ov{a}x_1 \cr
                            x_2\ov{a} & \ov{x}_1a & \xi_3}.$$
Then, $\delta_a \in F_4$.
\vspace{1mm}

(4) For $T \in O(3) = \{T \in M(3, \R) \, | \, {}^tTT = E \}$, we define a transformation $\delta(T)$ of $\gJ(3, \gC)$ by
$$
         \delta(T)X = TXT^{-1}, \quad X \in \gJ(3, \gC), $$
then, $\delta(T) \in F_4$.
\vspace{2mm}

Using (1) - (4), we will give a constructive proof of Theorem 1. Let 
$X = \pmatrix{\xi_1 & x_3 & \ov{x}_2 \cr
                        \ov{x}_3 & \xi_2 & x_1 \cr
                        x_2 & \ov{x}_1 & \xi_3}$  be a given element of $\gJ(3, \gC)$. 
\vspace{1mm}

(i) Assume $x_1 \not= 0$ and let $a = x_1/|x_1|$. Applying $\delta_a$ on $X$, then, the $x_1$-part of $X' = \delta_aX$ becomes real.
\vspace{1mm}

(ii) Applying some $T = \pmatrix{1 & 0 \cr
                                 0 & T'}, T' \in O(2) = \{T' \in M(2, \R) \, | \, {}^tT'T' = E \}$, then $X'$ can be transformed to the form $X'' = \pmatrix{{\xi_1}' & {x_3}'' & {\ov{x}_2}'' \cr
                        {\ov{x}_3}'' & {\xi_2}'' & 0 \cr
                        {x_2}'' & 0 & {\xi_3}''}$.
\vspace{1mm}

(iii) Assume ${x_2}'' \not= 0$ and let $a = {x_2}''/|{x_2}''|$. Applying $\delta_a$ on $X$, then, the $x_2$-part of $X^{(3)} = \delta_aX''$ becomes real. That is, $X^{(3)}$ is of the form $X^{(3)} = \pmatrix{{\xi_1}' & {x_3}^{(3)} & r_2 \cr                             {\ov{x}_3}^{(3)} & {\xi_2}'' & 0 \cr
                                  r_2 & 0 & {\xi_3}''},$ $ r_2$ is real.
                                                
\vspace{1mm}

(iv) Let $\C = \{x + ye_1 \, | \, x, y \in \R \}$ be the field of the complex numbers contained in $\H$ : $\C \subset \H \subset \gC$. Since the group $G_2$ acts transitively on ${S_r}^6 = \{u \in \gC \, | \, u = -\ov{u}, |u| = r \}$,  any element $x = x_0 + u \in \gC, x_0 \in \R, u \in \gC, u = -\ov{u}$ can be deformed to $c_3 = x_0 + x_1e_1 \in \C$ by some $\alpha \in G_2$ ([2]). Applying this $\alpha$ on $X^{(3)}$, then $\alpha X^{(3)}$ is of the form $X^{(4)} = \pmatrix{{\xi_1}' & c_3 & r_2 \cr
                                       \ov{c}_3 & {\xi_2}'' & 0 \cr
                                       r_2 & 0 & {\xi_3}''} \in \gJ(4, \C) \subset \gJ(3, \H)$.
 
(v) $X^{(4)}$ can be transformed to diagonal form by some $A \in Sp(3)$, that is, $AX^{(4)}A^*$ is of the form $\pmatrix{{\xi_1}'' & 0 & 0 \cr
                                 0 & {\xi_2}^{(3)} & 0 \cr
                                 0 & 0 & {\xi_3}^{(3)}}$. 
\vspace{1mm}

\noindent Note that all process of (i) - (v) are constructive. Thus, Theorem 1 is proved.
\vspace{3mm}

II. Let $\gC' = \ H \oplus \H{e_4}', {e_4}'^{\;2} = 1$ be the split Cayley algebra and $\gJ(3, \gC') = \{X \in M(n, \gC') \, | \, X^* $ $ = X \}$ with the Jordan multiplicaton 
$
      X \circ Y = \dfrac{1}{2}(XY + YX). $
Then, the connected non-compact Lie group $F_{4(4)}$ is defined by
$$
   F_{4(4)} = \{ \alpha \in \Iso_R(\gJ(3, \gC')) \, | \, \alpha(X \circ Y) = \alpha X \circ \alpha Y \}. $$
Then, we have the following Theorem.
\vspace{3mm}

{\bf Theorem 2.} {\it Any element $X$ of $\, \gJ'$ can not necessarily be transformed to a dinagonal form by element} $\alpha \in F_{4(4)}$.
\vspace{2mm}

{\bf Proof.} We will give a counter example. Assume that the element
$$
   X_0 = \pmatrix{0 & 0 & 0 \cr
                0 & 0 & {e_4}' \cr
                0 & -{e_4}' & 0} \in \gJ' $$ 
can be transformed to a diagonal form by some $\alpha \in F_{4(4)}$:
$$
   \alpha X_0 = \pmatrix{\xi_1 & 0 & 0 \cr
                         0 & \xi_2 & 0 \cr
                         0 & 0 & \xi_3}, \quad \xi_i \in \R. $$
If we define an inner product $(X, Y)$ in $\gJ(3, \gC')$ by $(X, Y) = \tr(X \circ Y)$, then we know that any element $\alpha \in F_{4(4)}$ leaves the inner product invariant : $(\alpha X, \alpha Y) = (X, Y)$. Now, we have 
$$
     (\alpha X_0, \alpha X_0) = (X_0, X_0) = 2{e_4}'(-{e_4}') = -2. $$
On the other hand, 
$$
   (\alpha X_0, \alpha X_0) = \Big(\pmatrix{\xi_1 & 0 & 0 \cr
                         0 & \xi_2 & 0 \cr
                         0 & 0 & \xi_3}, \pmatrix{\xi_1 & 0 & 0 \cr
                         0 & \xi_2 & 0 \cr
                         0 & 0 & \xi_3} \Big) = {\xi_1}^2 + {\xi_2}^2 + {\xi_3}^2 \geq 0,$$
which contradicts the above. Therefore, $X_0$ can not be trasformed to a diagonal form.
\vspace{3mm}

{\bf References}

[1] H. Freudenthal,  Oktaven, Ausnahmegruppen und Oktavengeometrie, Math. Inst. Rijks univ. te Utrecht, 1961.

[2] I. Yokota, Linear simple Lie groups of exceptional type (in Japanese), Gendai-Sugakusha, (1992).

[3] I. Yokota, Realizations of involutive automorphisms $\sigma$ and $G^\sigma$ of exceptional linear Lie groups $G$, Part I, $G_2, F_4$ and $E_6$, Tsukuba J. Math., 14(1990), 185 - 223.

\end{document}